\documentclass[10pt]{amsart}

\usepackage[all]{xy}
\usepackage{hyperref}

\newtheorem{proposition}{Proposition}[section]
\newtheorem{lemma}[proposition]{Lemma}
\newtheorem{corollary}[proposition]{Corollary}
\newtheorem{theorem}[proposition]{Theorem}

\theoremstyle{definition}
\newtheorem{definition}[proposition]{Definition}

\newtheorem{blanco}[proposition]{}

\theoremstyle{remark}
\newtheorem{remark}[proposition]{Remark}

\newcommand{\thlabel}[1]{\label{th:#1}}
\newcommand{\thref}[1]{Theorem~\ref{th:#1}}
\newcommand{\selabel}[1]{\label{se:#1}}
\newcommand{\seref}[1]{Section~\ref{se:#1}}
\newcommand{\lelabel}[1]{\label{le:#1}}
\newcommand{\leref}[1]{Lemma~\ref{le:#1}}
\newcommand{\prlabel}[1]{\label{pr:#1}}

\newcommand{\colabel}[1]{\label{co:#1}}
\newcommand{\coref}[1]{Corollary~\ref{co:#1}}
\newcommand{\relabel}[1]{\label{re:#1}}

\newcommand{\delabel}[1]{\label{de:#1}}
\newcommand{\deref}[1]{Definition~\ref{de:#1}}
\newcommand{\eqlabel}[1]{\label{eq:#1}}
\newcommand{\equref}[1]{(\ref{eq:#1})}

\newcommand{\ul}[1]{\underline{#1}}

\def\Z{\mathbb Z}

\def\R{\mathcal R}
\def\S{\mathcal S}
\def\*C{{^*\mathcal C}} 
\def\C{\mathcal C}
\def\D{\mathcal D}
\def\J{\mathcal J}
\def\M{\mathcal M}
\def\T{\mathcal T}

\newcommand{\Hom}{{\rm Hom}}
\newcommand{\End}{{\rm End}}

\newcommand{\Ker}{{\rm Ker}\,}

\newcommand{\im}{{\rm Im}\,}

\def\ot{\otimes}

\def\ZZ{{\mathbb Z}}

\newcommand{\Aa}{\mathcal{A}}

\newcommand{\Cc}{\mathcal{C}}
\newcommand{\Dd}{\mathcal{D}}
\newcommand{\Ee}{\mathcal{E}}

\newcommand{\Mm}{\mathcal{M}}

\newcommand{\Rr}{\mathcal{R}}
\def\*C{{}^*\hspace*{-1pt}{\Cc}}

\def\text#1{{\rm {\rm #1}}}

\def\ul{\underline}

\begin{document}
\title[Local units versus local projectivity]{Local units versus local projectivity. \\
Dualisations : Corings with local structure maps.}
\author{J. Vercruysse}
\address{Faculty of Applied Sciences,
Vrije Universiteit Brussel, VUB, B-1050 Brussels, Belgium}
\email{joost.vercruysse@vub.ac.be}
\urladdr{http://homepages.vub.ac.be/\~{}jvercruy/}
\thanks{}


\begin{abstract}
We unify and generalize different notions of local units and local projectivity. We investigate the connection between these properties by constructing elementary algebras from locally projective modules. Dual versions of these constructions
are discussed, leading to corings with local comultiplications, corings with local counits and rings with local multiplications. 
\end{abstract}

\maketitle

\section*{Introduction}
Corings and their comodules were introduced by Sweedler in \cite{Sweedler65}; there has been
a revived interest in the subject recently, after an observation made
by Takeuchi that Hopf modules and most of their generalizations, including for
example
Doi-Hopf modules, entwined modules, weak Hopf modules and Yetter-Drinfeld
modules are examples of comodules over certain corings. \cite{Br3} is the first of a series of papers with new applications
of corings; a detailed discussion appeared recently in \cite{BrzezinskiWisbauer}.

Rings without unit have been examined in \cite{Ab,AM}; various applications
have appeared in the literature. In \cite{Taylor}, a Brauer group of equivalence classes
of Azumaya algebras without unit was introduced (see also \cite{C}).
There is also a close connection between rings with local units and
the different notions of local projectivity that exist in the literature.
A first notion is due to Zimmermann-Huisgen \cite{ZH}, and is equivalent
to the so-called $\alpha$-condition (see for example \cite{W}).
It can be defined relative to a an abelian subgroup $\R$ of the dual
module ${^*M}={_A\Hom}(M,A)$. In the case where $M$ is an $A$-coring, there is a remarkable
relation between  $\R$-relative local projectivity of $\C$ and the
existence of local units in $\R$ (see \cite{CV}). A second notion (which we will call
strongly local projectivity) is due to Abrams \cite{Ab}, and is related to rings with idempotent local units by Morita theory (see \cite{AM}).

The aim of this paper is a further investigation of the relations between
local projectivity and local units. Let $B$ be a ring with unit, and
$M$ a right module over a $B$-ring $R$ without unit. We introduce the notion
of (idempotent) local unit map on $M$, and discuss how this generalizes the
local units of \cite{Ab,AM}. In \seref{3}, we look at the following situation:
we consider two rings $A$ and $B$ with units, and a dual pair of bimodules,
consisting of a $(B,A)$-bimodule $M$, an $(A,B)$-bimodule $M'$, and an $A$-bimodule
map $\mu:\ M'\ot_B M\to A$. We then discuss the relation between 
weak local projectivity of $M$ and the existence of local units on $M$
considered as a right module over the elementary algebra $M\ot_A M'$, which is
a $B$-ring (see Theorems \ref{th:3.4} and \ref{th:3.6}). We have similar
characterizations of strong local projectivity, but now in terms of
idempotent local units, and we have a structure theorem for dual pairs
of strongly locally projective modules (see \coref{3.8}).

Our next aim is to introduce similar notions for corings without counits
(but with a coassociative comultiplication), and to relate these to
properties of rings with local units. This leads us to corings with
(weak or strong) local comultiplications (\seref{4}) and corings with (idempotent) local counits
(\seref{5}). Both notions are related to local units:
if a coring has (weak or strong) local comultiplications, then
its endomorphism ring has (idempotent) local units on the coring considered
as a module over the endomorphism ring (\thref{4.2}); right counits on a
comodule $M$ are precisely the units over $M$ viewed as a module over
the dual ring (\thref{5.3}). Finally, there exists a dual notion of local comultiplication:
in \seref{7}, we introduce local multiplications. 

If $X$ is an object of a category, then $X$ will be also a notation for the
identity isomorphism on $X$.

\section{Preliminary results}\selabel{1}
\subsubsection*{Rings and corings}
Let $A$ be a ring, not necessary with unit. A ring with unit will be
called a unital ring.
$\M_A$ will be the category 
of right $A$-modules and $A$-module morphisms. We use similar notation
for the categories of left $A$-modules and $A$-bimodules.

An $A$-ring is an object $B\in{}_A\M_A$ together with an $A$-bimodule
map $\mu:\ B\otimes_A B \to B$ satisfying the associativity condition
 $\mu\circ (\mu\otimes_A I_B)=\mu
\circ (I_B\otimes_A \mu)$. This makes $B$ into a ring with associative
multiplication $\mu$. If $e\in B$ is a unit for this multiplication, then
the map $\iota:\ A\to B$, $\iota(a)=ae=eae=ea$ is a morphism of rings,
and a morphism of unital rings
if $A$ has a unit. By ring we mean a $\ZZ$-ring.

A module over the $A$-ring $B$, or shortly a $B$-module, is 
$M\in\M_A$ together with a right $A$-module map $\mu_M:\ M\otimes_AB\to 
M$ satisfying the usual associativity condition.
Observe that a module $M$ over the ring $B$ is not always a module over the 
$A$-ring $B$, since $M$ has no canonical $A$-module structure 
if $B$ has no unit map.

A right $B$-module is called \emph{firm} if the map 
$M\otimes_BB\to M$, $m\otimes_Bb\mapsto m\cdot b$ is a right $B$-module 
isomorphism.

Let $\Cc$ be an $A$-bimodule. A \emph{comultiplication} on $\C$ is an
$A$-bimodule map $\Delta_\C:\ \C\to\C\otimes_A\C, c\mapsto 
c_{(1)}\otimes_Ac_{(2)}$ such that this comultiplication is coassociative, 
that is,
$$\Delta(c_{(1)})\otimes_Ac_{(2)}=c_{(1)}\otimes_A\Delta(c_{(2)})
:=c_{(1)}\otimes_Ac_{(2)}\otimes_Ac_{(3)},$$
for all $c\in \C$.
An $A$-coring $\C$ is an 
$A$-bimodule together with a comultiplication and an $A$-bilinear map 
$\varepsilon_\C:\C\to A$, such that 
$$c_{(1)}\varepsilon_\C(c_{(2)})=\varepsilon_\C(c_{(1)})c_{(2)}=c,$$
for all $c\in\C$. We call $\varepsilon_\C$ a counit.

Let $\C$ and $\D$ be two $A$-corings, then an $A$-bilinear map $f:\C\to 
\D$ is called a morphism of $A$-corings if 
$\varepsilon_\D(f(c))=\varepsilon_\C(c)$ and 
$\Delta_\D(f(c))=f(c_{(1)})\otimes_A f(c_{(2)})$, for all $c\in \C$.

The left dual $\*C={_A\Hom}(\C,A)$ of an $A$-coring 
is an $A$-ring with multiplication
$$(f* g)(c)=g(c_{(1)}f(c_{(2)}))$$
and unit $\varepsilon_\C$.

A right $\C$-comodule $M$ is a right $A$-module, together with a right 
$A$-module map $\rho:\ M\to M\otimes_A\C$, $\rho(m)=m_{[0]}\otimes_Am_{[1]}$, 
such that 
$\rho(m_{[0]})\otimes_Am_{[1]}=m_{[0]}\otimes_A\Delta(m_{[1]}):=
m_{[0]}\otimes_Am_{[1]}\otimes_Am_{[2]}$
for all $m\in M$. We call $M$ counital if, in addition, $m=m_{[0]}\varepsilon_\C(m_{[1]})$ for all $m\in M$. A map 
$f:\ M\to N$ between two 
right $\C$-comodules is a $\C$-comodule morphism if it is a right 
$A$-linear, and right $\C$-colinear, by which we mean that
$$f(m)_{[0]}\otimes_Af(m)_{[1]}=f(m_{[0]})\otimes_Am_{[1]},$$
for all $m\in M$. $\M^\C$ is the category of right $\C$-comodules 
and $\Cc$-comodule morphisms. In a similar way, we define the categories
${^\C\M}$, ${^\C\M^\C}$, ${_A\M^\C}$, etc.

Let $\Cc$ be an $A$-coring, not necessarily with a counit, and
consider $M\in\Mm^\Cc$ and $N\in{^\Cc\Mm}$. The cotensor product
$M\otimes^\Cc N$ is the equalizer of $\rho_M\otimes_A I_N$ and $I_M\otimes_A \rho_N$.

\subsubsection*{Split direct systems}
Let $(C_i)_{i\in I}$ be a direct system in an Abelian category $\Aa$.
Then $I$ is a partially ordered set, such that for all $i,j\in I$, there exists
$k\in I$ such that $k\geq i$ and $k\geq j$ (we will write
$k\geq i,j$), and we have morphisms
$\varphi_{ji}:\ C_i\to C_j$ for $i\le j$, 
such that $\varphi_{ii}={C_i}$ and $\varphi_{kj}\circ\varphi_{ji}=\varphi_{ki}$ if $i\le j\le k$. $(C_i)_{i\in I}$ is called a \emph{split direct system} if
every $\varphi_{ji}$ has a left inverse $\psi_{ij}$, such that
$\psi_{ij}\circ\psi_{jk}=\psi_{ik}$ if $i\le j\le k$. It then follows that
$\psi_{ii}={C_i}$, and the $\varphi_{ji}$ are cosplit monomorphisms, and the
$\psi_{ij}$ are projections.

\begin{proposition}\prlabel{1.1}
Let $(C_i)_{i\in I}$ be a direct system in $\Aa$, and consider the direct limit
$$C=\lim_{\longrightarrow} C_i,$$
and the canonical maps $\varphi_i:\ C_i\to C$.
The direct system is split if and only if there exists for every $i\in I$
a morphism $\psi_i:\ C\to C_i$ such that 
\begin{eqnarray}
\psi_i\circ \varphi_i&=&C_i; \\
\psi_i\circ \varphi_j\circ\psi_j&=&\psi_i, \qquad \textrm{for every } i\leq j.
\eqlabel{psi}
\end{eqnarray}
\end{proposition}

\begin{proof}
For our purposes, it is sufficient to consider the case $\Aa=\Mm_A$.
Recall first (see for example \cite{St}) that $C$ is the disjoint
union of the $C_i$, modulo the equivalence relation $\sim$ defined as follows:
for $c_i\in C_i$ and $c_j\in C_j$, $c_i\sim c_j$ if and only if there exists
$k\geq i,j$ such that $\varphi_{ki}(c_i)=\varphi_{kj}(c_j)$. 
$\varphi_i:\ C_i\to C$ is then given by $\varphi_i(c)=[c]$, with $c\in C_i$ and $[c]$ the equivalence class of $c$.

Assume now that we have a split direct system. We define $\psi_i:\ C\to C_i$
as follows: assume that $x\in C$ is represented by $c_j\in C_j$.
Then take $k\geq i,j$, and put
$$\psi_i(x)=\psi_i\circ\varphi_j(c_j):=(\psi_{ik}\circ\varphi_{kj})(c_j).$$
Let us show that $\psi_i$ is well-defined. First we show that the definition
is independent of the choice of $k\geq i,j$. Take $l\geq i,j$ and then
$m\geq i,j$. We then have
$$
(\psi_{ik}\circ\varphi_{kj})(c_j)=(\psi_{ik}\circ\psi_{km}\circ \varphi_{mk}\circ \varphi_{kj})(c_j)=(\psi_{im}\circ\varphi_{mj})(c_j),$$
and, in a similar way
$$(\psi_{il}\circ\varphi_{lj})(c_j)=(\psi_{im}\circ\varphi_{mj})(c_j).$$
Now we show that the definition of $\psi_i$ is independent of the choice
of the representing $c_i$. Let $x=[c_j]=[c_l]$, with $c_j\in C_j$ and
$c_l\in C_l$. Then take $k\geq i,j,l$ such that $\varphi_{kj}(c_j)=
\varphi_{kl}(c_l)$. Then $(\psi_{ik}\circ \varphi_{kj})(c_j)=
(\psi_{ik}\circ \varphi_{kl})(c_l)$, as needed.\\
Now take $c_i\in C_i$. Then
$$\psi_i(\varphi_i(c_i))=\psi_{ii}(c_i)=c_i.$$
To prove the second equality, we first check that for every $i\leq j$,
$$\psi_{ij}=\psi_{ij}\circ C_j=\psi_{ij}\circ\psi_{jk}\circ\varphi_{kj}
=\psi_{ik}\circ\varphi_{kj}=\psi_i\circ\varphi_j,$$
where we took $k\geq i,j$.
Let $l\leq i$, and take $x=[c_j]\in C$, $k\geq i,j,l$. Then we compute
$$(\psi_l\circ\varphi_i\circ\psi_i)(x)=(\psi_{li}\circ \psi_i)(x)=(\psi_{li}\circ \psi_{ik}\circ \varphi_{kj})(c_j)
=(\psi_{lk}\circ \varphi_{kj})(c_j)=\psi_l(x).$$
Conversely, define $\psi_{ij}=\psi_i\circ\varphi_j$. Then we find
$$\psi_{ij}\circ\psi_{jk}=\psi_i\circ\varphi_j\circ\psi_j\circ\varphi_k
=\psi_i\circ\varphi_k=\psi_{ik},$$
and
$$\psi_{ij}\circ\varphi_{ji}=\psi_i\circ\varphi_j\circ\varphi_{ji}=
\psi_i\circ\varphi_i={C_i}.$$
\end{proof}

\begin{remark}
From the proof it is clear that under the equivalent conditions of the previous proposition, \equref{psi} means exactly
$$\psi_i=\psi_{ij}\circ\psi_j,$$
for every $i\leq j$.
\end{remark}

\section{Rings with local units}\selabel{2}
\subsection{local units}\selabel{2.1}

\begin{definition}\delabel{2.1}
Let $B$ be a ring with unit and $R$ a $B$-ring. Let $M$ be a right
module over the $B$-ring $R$ (this implies $M$ is a right 
$B$-module by definition).
A \emph{right unit map} on $M$ is a $B$-bimodule map
$\eta_M:\ B\to R$, 
such that the following diagram of right $B$-modules is commutative 
\[
\xymatrix{
M\otimes_B B \ar[d]_{I_M\otimes_B\eta_M} \ar[dr]^{\cong} \\
M\otimes_B R \ar[r]_{~~~\mu_M} & M
}
\]
Left unit maps are defined in a similar way. If 
$M$ is a bimodule and $\eta_M$ a left and right unit map, then we 
call $\eta_M$ just a unit map.
\end{definition}

\begin{lemma}\lelabel{2.2}
Let $R$ be a $B$-ring and $M$ a right $R$-module.
The following statements are equivalent.
\begin{enumerate}
\item For every $m\in M$ there exists an element $e\in R^B:=\{r\in 
R~|~br=rb,~{\rm for~all}~ b\in B\}$ such that $m\cdot e=m$;
\item for every finitely generated $B$-submodule $T$ of $M$, there exists 
an element $e\in R^B$ such that $t\cdot e=t$ for all $t\in T$;
\item there exists a 
unit map $\eta_T$ on every finitely generated $B$-submodule $T$ of $M$.
\end{enumerate}
In this case we say that $R$ has right local 
units on $M$, and we call $e$ a right local unit and $\eta_T$ a right 
local unit map on $T$.
\end{lemma}

\begin{proof}
$\underline{1\Rightarrow 2}$ Let $\{t_1,\ldots,t_k\}$ be a set of 
generators for $T$. We proceed by induction on $k$. If $k=1$ then $T=tB$ 
and from the first statement, we find an $e\in R^B$ such that $te=t$. 
Consequently $tbe=teb=tb$. If $k>1$, we can find a right local unit $e'$ 
for $t_1$. From the inducion hypothesis, we can also find a right local 
unit $e''$ for the $k-1$ elements $t_i-t_ie'$, $i=2,\ldots,k$. Now 
$e:=e'+e''-e'e''\in R^B$ and
\begin{eqnarray*}
t_1e&=&t_1e'+t_1e''-t_1e'e''=
t_1+t_1e''-t_1e''=t_1;\\
t_ie&=&t_ie'+t_ie''-t_ie'e''=
t_ie'+(t_i-t_ie')e''=
t_ie'+t_i-t_ie'=t_i,
\end{eqnarray*}
and we find that $e$ is a right local unit for all $t_i$ ($1\le i\le k$) 
and consequently for every element in $T$.\\
$\underline{2\Rightarrow 1}$ is trivial.\\
$\underline{2\Leftrightarrow 3}$ follows from the fact 
that $R^B\cong {_B\Hom_B}(B,R)$ as $\Z$-modules. 
\end{proof}

\begin{definition}\delabel{2.3}
A $B$-ring $R$ has right \emph{local units}, if it has right local units as
a right module over itself. In the same way, we define rings with left and two-sided local units. By local units, we will always mean two-sided local units.\\
If $R$ has right local units, then a right $R$-module $M$ on which
$R$ has right local units will be termed a right \emph{unital} $R$-module.
\end{definition}

\begin{lemma}\lelabel{2.5}
Let $R$ be a $B$-ring, $M$ a right $R$-module and $N$ a left $R$-module. 
If there exist a right local 
unit $e'\in R$ on $m\in M$ and a left local unit $e''\in R$ on $n\in N$, 
then there exists also an element $e\in R$ which is at the same time a 
right local unit on $m\in M$ and a left local unit on $n\in N$.
\end{lemma}

\begin{proof} Take $e:=e'+e''-e'e''$.
\end{proof}

\begin{corollary}\colabel{2.6}
If $R$ has left and right local units on a bimodule $M$, then it also has 
two-sided local units on $M$. In particular, if $R$ is a $B$-ring with left 
and right local units, then $R$ is a $B$-ring with two-sided local units.
\end{corollary}

\begin{remark}\relabel{2.4}
Recall from \cite{T} that a ring $R$ is called left (resp. right) $s$-unital
if $u\in Ru$ (resp. $\in uR$), for every $u\in R$. $R$ is called $s$-unital if $u\in Ru\cap uR$, for every $u\in R$. It follows immediately
from Lemmas \ref{le:2.2} and \ref{le:2.5} and \coref{2.6} that a $\ZZ$-ring $R$
has left (resp. right, two-sided) local units if and only if it is
left $s$-unital (resp. right $s$-unital, $s$-unital).
\end{remark}

\begin{lemma}\lelabel{2.7}
Let $R$ be a $B$-ring with right local units. A right $R$-module $M$ has
right local units if and only if it is firm.
\end{lemma}

\begin{proof}
Assume that $M$ has right local units; the map
$\psi:\ M\to M\otimes_R R$, $\psi(m):= m\otimes_R e$, with $e$ a 
right local unit on $m$, is well-defined. Let $e'$ be another 
right local 
unit on $m$, and choose a right local unit $e''$ on $e$ and $e'$. We find 
$m\otimes_R e=m\otimes_R ee'' = me\otimes_R e'' = m\otimes_R 
e''=me'\otimes_R e'' = m\otimes_R e'e'' = m\otimes_R e'$. It is obvious
that $\psi$ is inverse to the module structure map $M\otimes_R R\to M$.

Conversely, if $M$ is a firm right
$R$-module, then $M\cong M\otimes_R R$, and since $R$ has local units on 
itself, it also has local units on $M$.
\end{proof}

\begin{theorem}\thlabel{2.8}
\label{locmod}
If $R$ is a $B$-ring with right local units,
then every firm right module $M$ over $R$ (considered as a $\ZZ$-ring) is 
also a right $B$-module. Moreover the category of firm right modules over 
$R$ as a $\ZZ$-ring is isomorphic to the category of firm right modules over $R$ as 
$B$-ring. 
\end{theorem}

\begin{proof}
We know that $M\cong M\otimes_R R$, and since $R$
is a right $B$-module, $M$ is a right $B$-module as well. The 
$B$-action is given by the formula $m\cdot b= m(eb)$, with $e$ a
right local unit on $m$.
A similar argument applies to the morphisms.
\end{proof}

\subsection{Idempotent local units}\selabel{2.2}
\begin{definition}
With notation as in \deref{2.1}, a right local unit map $\eta_M$
on $M$ is called \emph{idempotent} if $\eta_M$ is a morphism of $B$-rings.
This means that the following diagram commutes:
\[
\xymatrix{
B\otimes_B B \ar[rr]^{\eta_M\otimes_B\eta_M} \ar[d]_{\mu_B} 
& & R\otimes_B R \ar[d]^{\mu_R}\\
B \ar[rr]_{\eta_M} & & R
}
\]
\end{definition}

\begin{lemma}
\lelabel{2.10} Let $R$ be a $B$-ring, and $M$ a right $R$-module. 
The following statements are equivalent.
\begin{enumerate}
\item $R$ has an idempotent right local unit map on every finitely generated right $B$-submodule $N$ of $M$;
\item for every finitely generated right $B$-submodule $N$ of $M$, we can 
find an idempotent $e\in R^B$ such that $ne=n$ for all $n\in N$;
\item for every finitely generated right $B$-submodule $N$ of $M$, we can 
find an idempotent $e\in R^B$ such that $N\subseteq Me$;
\item $M$ is the direct limit of a split direct system of submodules $(M_i)_{i\in I}$,
and $R$ has an idempotent right unit map on every $M_i$.
\end{enumerate}
In this case we say that $R$ has idempotent 
right local units on $M$ and we call $e$ an idempotent right local unit.
\end{lemma}

\begin{proof}
The equivalence $1\Leftrightarrow 2$ can be proved in the same way as \leref{2.2}, 
taking into account that ring morphisms $B\to R$ 
correspond to idempotents in $R^B$. 

$\ul{3\Leftrightarrow 2}$ is trivial.

$\ul{4\Rightarrow 2}.$ Let $\{n_1,\ldots,n_k\}$ be a 
set of generators of a finitely generated right $B$-submodule $N$
of $M$. Since $M=\varinjlim M_i$, 
there exists $\ell\in I$ such that $n_i\in M_l$, for all $j=1,\cdots,k$.
The idempotent right unit map on $M_\ell$ is an idempotent right local unit map on $N$.

$\ul{2\Rightarrow 4}.$ Let $I$ be the set of all idempotent elements $e\in R^B$.We define a partial order as follows: $e\le f$ if $ef = fe = e$. Setting $M_e:=Me$, then $\{M_e\}_{e\in I}$ is
a split direct system, with maps
$$\varphi_{fe}:\ M_e\to M_f,~~\varphi_{fe}(me)=me=mef;$$
$$\psi_{ef}:\ M_f\to M_e,~~\psi_{ef}(mf)=mfe=me,$$
for all $m\in M$, and $e\leq f\in I$. $e$ is an idempotent local unit on $Me$,
and $M=\varinjlim M_e$.
\end{proof}

It follows from \leref{2.10} that $R$ is a $B$-ring with idempotent left
(resp. right) local units if and only if $R$ is the direct limit of a split
direct system of rings with a left (resp. right) unit.

\begin{lemma}\lelabel{2.11}
A $B$-ring $R$ with left and right idempotent local
units also has two-sided idempotent local units. 
\end{lemma}

\begin{proof}
For every finite subset $\{r_1,\ldots,r_k\}\subset R$, we have to find an 
idempotent $e\in R^B$ such 
that $er_i=r_ie=r_i$. By assumption, we can find an idempotent left 
local unit $e'$ on $\{r_1,\ldots,r_k\}$, and an idempotent right
local unit $e''$ on $\{r_1,\ldots,r_k,e'\}$. An easy calculation shows that
$e=e'+e''-e''e'$ is an idempotent two-sided unit on $\{r_1,\ldots,r_k\}$.
\end{proof}

\begin{remark}
It follows from Lemmas \ref{le:2.10} and \ref{le:2.11} that
$R$ is a $B$-ring with 
two-sided local units if and only if for every finitely generated $B$-subbimodule
$F$ of 
$R$, there exists an idempotent $e\in R^B$ such that $F\subseteq eRe$. Notice that 
$eRe$ is a $B$-ring with unit $e$.
In the case where $B=\Z$, we recover the definition of ring with local units as 
introduced by \'Anh and M\'arki in \cite{AM}. 
\end{remark}

\subsection{Local Projectivity}\selabel{2.3}
In order to be able to distinguish different notions of local projectivity, we
introduce ``weak" and ``strong" local projectivity. In the literature,
both notions are termed ``local projectivity".

Let $A$ and $B$ be a rings with unit and $M\in {}_B\M_A$. Then 
$M^*=\Hom_A(M,A)$ is an $(A,B)$-bimodule, with actions 
$(afb)(m)=af(bm)$ for all $f\in M^*$, $m\in M$, $a\in A$ and $b\in B$.

\begin{definition}\delabel{2.13}
With notation as above, let $\Rr$ be an additive subgroup of $M^*$, 
and $N$ a subset of $M$. A \emph{dual ($M,\R$)-basis} of $N$ is a finite 
set $\{(u_i,f_i)\}\subset M\times \R$ such that $n=\sum_i u_if_i(n)$, for 
every $n\in N$. 

$\{(u_i,f_i)\}$ is called (left) $B$-linear if $\sum_i u_if_i(bm)=\sum_i bu_if_i(m)$, for all $m\in M$.

We call $\{(u_i,f_i)\}$ {idempotent} if $u_i=\sum_j u_jf_j(u_i)$.

We call $M$ \emph{weakly $\Rr$-locally projective} as a $(B,A)$-bimodule if every 
finite subset $N$ of $M$ has a $B$-linear dual 
$(M,\R)$-basis. If $\Rr=M^*$ then we just say that $M$ is right weakly locally projective as a $(B,A)$-bimodule.
\end{definition}

\begin{lemma}\lelabel{2.14}
Let $M$ be a $(B,A)$-bimodule.
$M$ is weakly $\Rr$-locally projective as $(B,A)$-bimodule if
and only if every finitely generated $(B,A)$-subbimodule of $M$ has a $B$-linear dual 
$(M,\R)$-basis. 
\end{lemma}

\begin{proof} This follows by the $(B,A)$-linearity of the dual bases, where the $A$-linearity is satisfied since $\Rr\subset \Hom_A(M,A)$.
\end{proof}

In case $B=\Z$, we say that $M$ is weakly $\Rr$-locally projective as a right
$A$-module. From \cite{CV}, we recall the following characterization, which has
an obvious analog for $(B,A)$-bimodules.

\begin{theorem}\thlabel{2.15}
With notation as above,  the following statements are equivalent 
\begin{enumerate}
\item $M$ is weakly $\R$-locally projective as a right $A$-module;
\item $M$ is weakly $\S$-locally projective as a right $A$-module, where $\S$ 
is the left $A$-submodule of $M^*$ generated by $\Rr$;
\item the map $\alpha_{N,\R}:\ M\otimes_AN\to {}_A\Hom(\R,N)$,
$\alpha(m\otimes_A n)(f)=f(m)n$, is injective for every left $A$-module $N$, that is, $M$ satisfies the $\alpha$-condition for $\R$;
\item $M$ is weakly $\T$-locally projective and $\S$ is dense in the finite topology on $M^*$ for every $A$-submodule $\T$ of $M^*$ such that $\S\subset\T$.
\end{enumerate}
\end{theorem}

$M$ is weakly $M^*$-locally projective as a
right $A$-module if and only if $M$ is locally projective in the sense of
Zimmermann-Huisgen \cite{ZH}. This is equivalent to the following condition
(see \cite{G}): for any commutative diagram with exact rows
in the category of right $A$-modules of the form
\[
\xymatrix{
0 \ar[r] & F \ar[r]^i & M \ar[d]^g \\
& N' \ar[r]_f & N \ar[r] & 0
}
\]
with $F$ finitely generated, there exists a right $A$-linear map 
$h:\ M\to N'$ such that $g\circ i = f\circ h\circ i$.

It is clear that projective modules are weakly locally projective and 
that weakly $\R$-locally projective modules are weakly locally 
projective.

\begin{definition}\delabel{2.16}
We call $M$ \emph{strongly $\R$-locally projective} as a $(B,A)$-bimodule if every 
finite subset $N$ of $M$ has an idempotent $B$-linear $(M,\R)$-dual basis. 

If $B=\Z$, we say that $M$ is strongly $\R$-locally projective as a right 
$A$-module. If $\R=M^*$, we say that $M$ is strongly locally 
projective a $(B,A)$-bimodule.
\end{definition}

It is clear that $M$ is strongly $\R$-locally projective as a $(B,A)$-bimodule if every 
finitely generated $(A,B)$-subbimodule $N$ of $M$ has an idempotent $B$-linear 
$(M,\R)$-dual basis.

It is clear that strongly locally projective modules are weakly locally 
projective. The converse implication is not true in general, since projective
modules are not necessarily strongly locally projective.

In \cite{AM}, a right $A$-module $P$ is called locally projective if it is
the direct limit of a split direct system $(P_i)_{i\in I}$ consisting of
submodules that are finitely generated and projective as right $A$-modules.

\begin{theorem}\thlabel{2.17}
The following statements are equivalent:
\begin{enumerate}
\item $P$ is locally projective in the sense of \'Anh and M\'arki \cite{AM};
\item every finitely generated submodule $F$ of $P$ is contained in a 
finitely generated projective submodule $P_F$ of $P$, which is 
a direct summand of $P$;
\item $P$ is strongly locally projective as a right $A$-module.
\end{enumerate}
\end{theorem}

\begin{proof}
$\underline{1\Rightarrow 2}$. Let $\{p_1,\ldots,p_n\}$ be a 
set of generators of $F$. Since $P=\varinjlim P_i$, there exists
$k\in I$ such that $\{p_1,\ldots,p_n\}\subset P_k$, hence
$F\subset P_k$, and $P_k$ is finitely generated projective and a
direct summand of $P$ as a right $A$-module.\\

$\underline{2\Rightarrow 3}$. Let $\{(u_k,f_k)\}$ be a dual basis of $P_F$
as a right $A$-module. Since $P_F$ is a direct summand of $P$, the 
maps $f_k:\ P_F\to A$ can be extended to $f_k:\ P\to A$.\\

$\underline{3\Rightarrow 1}$. Let $\{p_i\}_{i\in I}$ be a set of
generators of $P$ and let $\J$ be the set of all finite subsets of $I$.
For every $J\in \J$, let $\{(u^J_k,f^J_k)\}$ be a finite dual basis of the
module generated by $\{p_j\}_{j\in J}$. We define $P_J$ as the finitely
generated and projective submodule of $P$ generated by the $\{u_k^J\}$ and
$\psi_J:P\to P_J, \psi(p)=\sum u_kf_k(p)$.  Clearly these are projections and
consequently the $P_J$ are direct summands of $P$.  Furthermore, we define $J\le J'$ iff $\{(u^{J'}_{k'},f^{J'}_{k'})\}$ is a dual basis for the elements 
$u^J_k$. Since the dual bases are idempotent we have $J\le J$ and everything else is now straightforward.
\end{proof}

\section{Local projectivity versus local units}\selabel{3}
Let $A$ and $B$ be rings with unit. Following \cite[Sec. 1.2]{C},
a {dual pair} is a triple 
$\ul{M}=(M,M',\mu)$, with
$M\in {}_B\M_A$, $M'\in {}_A\M_B$ and $\mu:\ M'\otimes_B M\to A$
is an $A$-bimodule map.

Recall that $M^*=\Hom_A(M,A)$ is an $(A,B)$-bimodule, and that
${}^*M'={}_A\Hom(M',A)$ is a $(B,A)$-bimodule. Then
$(M,M^*,\mu)$, with $\mu(f\otimes_B
m)=f(m)$ is a dual pair. Also recall that the adjunction property
gives us the following isomorphisms of $\ZZ$-modules:
$${}_A\Hom_A(M'\otimes_BM,A)\stackrel{\zeta}{\longrightarrow} 
{}_A\Hom_B(M',M^*)\stackrel{\xi}{\longrightarrow} 
{}_B\Hom_A(M,{^*M'}).$$
For later use, we give the explicit description of the connecting maps:
$$\zeta(\mu)(m')(m)=\mu(m'\ot_B m),~~\zeta^{-1}(\varphi)(m'\ot_B m)=
\varphi(m')(m),$$
for all $\mu\in {}_A\Hom_A(M'\otimes_BM,A)$,
$\varphi\in {}_A\Hom_B(M',M^*)$, $m\in M$ and $m'\in M'$. We can consider
$\varphi^*\in {}_B\Hom_A({^*(M^*)},{^*M'})$, and
$\xi(\varphi)=\varphi^*\circ \iota$, with $\iota :\ M\to {^*(M^*)}$
the canonical morphisms. The proof of \leref{3.2} is straightforward.

\begin{lemma}\lelabel{3.2}
Let $A$ and $B$ be rings with unit, and $\ul{M}=(M,M',\mu)$ a dual pair.
Then $S:=M\ot_AM'$ is a $B$-ring, with multiplicatin map
$M\ot_A \mu\ot_A M$, and we have morphisms of $B$-rings
$\Phi:\ S\to \End_A(M)$ and $\Psi:\ S\to {}_A\End(M')$, given by
$$\Phi(m\ot_A m')(n)=m\mu(m'\ot_A n)~~{\rm and}~~
\Psi(m\ot_A m')(n')=\mu(n'\ot_A m)m'.$$
This makes $M$ into a left $S$-module, and $M'$ into a right $S$-module.
$S$ is called the elementary $B$-ring associated to $\ul{M}$.
\end{lemma}

A finite dual basis is a set
$\{(u_i,f_i)~|~ 1\le i\le n\}\subset M\times M^*$, such that $m=
\sum_iu_if_i(m)$ for all $m\in M$. It can also be regarded as an
element $e=\sum_i u_i\otimes_A f_i\in
M\otimes_A M^*$, such that $\Phi(e)=M$, the identity on $M$. 

\begin{lemma}\lelabel{3.3}
If there exists $\sum_i u_i'\otimes_B u_i\in 
(M'\otimes_B M)^A$ such that $\mu(\sum_i u_i'\otimes_B u_i)=1_A$,
then $M$ is a firm left $S$-module, and
$M'$ is a firm right $S$-module.
\end{lemma}

\begin{proof}
For every
$m\in M$, we have that $m=m1_A=m\mu(u_i'\otimes_B u_i)=\Phi(m\otimes_Au_i')(u_i)=(m\otimes_Au_i')\cdot u_i$, and this shows
that the canonical map $S\ot_S M\to M$ is surjective. We still have to show
that this map is injective. To this end, it suffices to show that
the sequence
$$S\ot_B S\ot_B M \stackrel{\lambda}{\longrightarrow}S\ot_B M
\stackrel{\psi}{\longrightarrow} M\to 0,$$
with $\lambda(s\ot_B t\ot_B m)=st\ot_B m-s\ot_B tm$ and $\psi(t\ot_B m)=tm$,
is exact. It is clear that $\psi\circ \lambda=0$. Take
$$x=\sum_j m_j\ot_A m'_j\ot_B n_j\in \Ker \psi,$$
this means that $\sum_j m_j\mu(m'_j\ot_B n_j)=0$. Then
\begin{eqnarray*}
&&\hspace*{-20mm}
\lambda\Bigl( \sum_{i,j} m_j\ot_A u_i'\otimes_B u_i\ot_A m'_j\ot_B n_j\Bigr)\\
&=& \sum_j m_j\ot_A m'_j\ot_B n_j-
\sum_{i,j} m_j\ot_A u_i'\otimes_B u_i\mu(m'_j\ot_B n_j)\\
&=& \sum_j m_j\ot_A m'_j\ot_B n_j-
\sum_{i,j} m_j\ot_A \mu(m'_j\ot_B n_j)u_i'\otimes_B u_i\\
&=& \sum_j m_j\ot_A m'_j\ot_B n_j=x\in \im(\lambda).
\end{eqnarray*}
\end{proof}

Let $A$ and $B$ be rings with unit, and $M\in {}_B\Mm_A$, $M'\in {}_A\Mm_B$.
We then have functors $F=-\ot_B M$ and $G=-\ot_A M'$ between the categories
$\Mm_B$ and $\Mm_A$, and adjunctions $(F,G)$ correspond to comatrix coring
contexts $(A,B,M',M,\mu,\nu)$. These consist of data $A$, $B$, $M$, $N$ as
above, and two bimodule maps $\mu:\ M'\otimes_BM\to A$ and $\nu:\ B\to M\otimes_AM'$ such 
that ${M'}=(\mu\otimes_A{M'})\circ ({M'}\otimes_B\nu)$ and 
$M=(\nu\otimes_BM)\circ(M\otimes_A\mu)$. This observation (see for example
\cite[Theorem 1.1.3]{C}) is folklore; the terminology was introduced recently
in \cite{BGT}.

\begin{theorem}\thlabel{3.4}
Let $A$ and $B$ be rings with unit, and $\ul{M}=(M,M',\mu)$ be a dual
pair.
With notation as above, the following statements are equivalent
\begin{enumerate}
\item $S$ is a $B$-ring with unit and $M$ and $M'$ are firm,
respectively as a left and right $S$-module;
\item $M$ is finitely generated and projective as a right $A$-module and 
$\zeta(\mu)=\varphi:\ M'\to M^*$ is bijective;
\item $M'$ is finitely generated and projective as a left $A$-module and 
$\xi(\zeta(\mu))=\psi:\ M\to {}^*M'$ is bijective;
\item the map $\Phi$ of \leref{3.2} is an isomorphism;
\item the map $\Psi$ of \leref{3.2} is an isomorphism;
\item there exists a $B$-bimodule map $\eta:\ B\to M\ot_A M'$ such that
$(A,B,M',M,\mu,\eta)$ is a comatrix coring context.
\end{enumerate}
\end{theorem}

\begin{proof}
$\ul{1.\Rightarrow 2.}$ Let $e=\sum_i u_i\ot u'_i\in M\ot_A M'=S$
be the unit element of $S$. It follows from \leref{2.7} that $M$ and
$M'$ are unital, resp. as a left and right $S$-module. For all $m\in M$,
we have that
$$m=em=\sum_i u_i\mu(u'_i\otimes_B m)=\sum_i u_i\varphi(u'_i)(m),$$
so $\{(u_i,\varphi(u'_i))\}$ is a finite dual basis of $M$ as a right $A$-module.
Since $\varphi$ is left $A$-linear, we have for every $f\in M^*$ that 
$$f=\sum_if(u_i)\varphi(u'_i)=\varphi(\sum_i f(u_i)u'_i)\in \im(\varphi),$$ 
so $\varphi$ is surjective. If $\varphi(m')=0$, then
$\mu(m'\ot m)=0$, for all $m\in M$, hence
$m'=m'u=\mu(m'\ot u_i)u'_i=0$, and it follows that $\varphi$ is injective.\\
$\ul{2.\Rightarrow 4.}$ and $\ul{4.\Rightarrow 1.}$ are trivial.\\
The equivalence of 1., 3. and 4. can be proved in a similar way.\\
$\ul{1.\Rightarrow 6.}$ Take $\eta:\ B\to M\ot_A M'$, $\eta(b)=be$, with
$e$ the unit of $S$.\\
$\ul{6.\Rightarrow 1.}$ $\eta(1)$ is a unit on $S$, and acts trivially on
$M$ and $M'$.
\end{proof}

We will now discuss how properties of local units can be translated
into properties of local projectivity. A first result is the following.

\begin{theorem}\thlabel{3.6}
We keep the notation introduced above, and let $\R:=\im \zeta(\mu)$. The following 
statements are equivalent:
\begin{enumerate}
\item $M$ is a weakly (resp. strongly) $\R$-locally projective right $A$-module;
\item $S$ is a $\Z$-ring with left local units (resp. left idempotent 
local units) and $M$ is a firm left $S$-module;
\item $S$ has local units (resp. idempotent local units) on $M$ as $\Z$-ring.
\end{enumerate}
In this situation,
$\xi(\zeta(\mu))=\psi$ is injective.
\end{theorem}

\begin{proof}
$\ul{1.\Rightarrow 2.}$ We first show that $S$ has left local units on $M$:
take $F=\{m_1,\cdots,m_n\}\subset M$; then there exists a dual $(M,\Rr)$-basis 
$\{(u_i,\varphi_\mu(u_i'))\}$ of the submodule of $M$ generated by $F$.
Let $e=\sum_i u_i\ot_A u_i'\in S$; then we have,
for all $j\in \{1,\cdots, n\}$, that 
$e m_j=\sum_iu_i\varphi_\mu(u_i')(m_j)=m_j$.\\
A similar argument shows that $S$ has left local units: take
$G=\{s_1,\cdots,s_m\}\subset S$, and write $s_j=\sum_k m_{jk}\ot_A m'_{jk}$.
There exists a dual $(M,\Rr)$-basis 
$\{(u_i,\varphi_\mu(u_i'))\}$ of the submodule of $M$ generated by the
$m_{jk}$, and $\sum_i (u_i\ot_A u_i')s_j=s_j$, for all $j$.\\
It then follows from \leref{2.7} that $M$is firm as a left $S$-module.\\
If $M$ is strongly $\R$-locally projective, then we can choose an idempotent
dual $(M,\Rr)$-basis $\{(u_i,\varphi_\mu(u_i'))\}$ in the above construction.
Then
\begin{eqnarray*}
e^2&=&(\sum_i u_i\otimes_A 
u_i')(\sum_ju_j\otimes_A u_j')=\sum_{i,j}u_i\mu(u_i'\otimes_B u_j)\otimes_A 
u_j'\\
&=&\sum_{i,j}u_i\varphi_\mu(u_i')(u_j)\otimes_A u_j'=
\sum_ju_j\otimes_A u_j'=e,
\end{eqnarray*}
and the local units are idempotent.

$\ul{2.\Rightarrow 3.}$ follows from \leref{2.7}.

$\ul{3.\Rightarrow 1.}$ Take $m\in M$, and a local unit 
$e=u_i\otimes_A u_i'\in S$ on the submodule of $M$ generated by $m$. Then 
$$m=e\cdot m=\sum_iu_i\mu(u_i'\otimes_B 
m)=\sum_iu_i\varphi_\mu(u_i')(m),$$
 so 
$\{(u_i,\varphi_\mu(u_i'))\}$ is a dual basis for $m$.

Finally, if $m\in \Ker \psi$, then we have for all $m\in M$ that
$\psi(m)(m')=\varphi(m')(m)=0$. Take a dual $(M,\Rr)$-basis 
$\{(u_i,\varphi_\mu(u_i'))\}$ of $m$. then 
$m=u_i\varphi(u_i')(m)=0$, and it follows that $\psi$ is injective.
\end{proof}
 
\begin{theorem}\thlabel{3.7} 
We keep the notation from \thref{3.6}. 
If $\Phi$ is injective, then the following statements are equivalent:
\begin{enumerate}
\item $M$ is a weakly (resp. strongly) $\R$-locally projective 
$(B,A)$-bimodule;
\item $S$ is a $B$-ring with left local units (resp. right 
idempotent local units) and $M$ is a firm left $S$-module;
\item $S$ has local units (resp. idempotent local units) on $M$ as a $B$-ring.
\end{enumerate}
If $\zeta(\mu)=\varphi$ is injective and $M$ is a weakly $\R$-locally projective 
$(B,A)$-bimodule, then
$\Psi$ is injective and $S$ is a $B$-subring of ${}_A\End(M')$.
\end{theorem}

\begin{proof}
$\ul{1.\Rightarrow 2.}$ 
We only have to show that the local units 
$e=\sum_i u_i\ot_A u_i'$ constructed in the proof of $\ul{1.\Rightarrow 2.}$ 
in \thref{3.6} can be taken in $S^B$. We can take a $B$-linear dual $(M,\Rr)$-basis
$\{(u_i,\varphi(u_i')\}$ of $F$. We then have for all $j$ and $b\in B$:
$$
 \Phi(eb)(m_j)=\Phi(e)(bm_j)=\sum_i u_i\varphi(u_i')(bm_j)
 = \sum_i bu_i\varphi(u_i')(m_j)=\Phi(be)(m_j),$$
so $\Phi(eb)=\Phi(be)$, and $e\in S^B$, since $\Phi$ is injective.

$\ul{2.\Rightarrow 3.}$ follows by \leref{2.7}.

$\ul{3.\Rightarrow 1.}$ The $B$-linearity of the dual basis is an immediate consequence of the fact that the local units commute with the elements of $B$.

If $\varphi$ is injective, then we can identify $M'$ and 
$\R$. Since $M$ is weakly $\R$-locally projective, it satisfies the 
$\alpha$-condition for $\R$. This means $M\otimes_A N\to {}_A\Hom(\R,N)$ 
is injective. Taking $N=M'$, we find that $\Psi$ is injective.
\end{proof}

\begin{corollary}\colabel{3.8}
Put $\S=\im(\psi)$.
The following statements are equivalent
\begin{enumerate}
\item $M$ is a weakly (resp. strongly) $\R$-locally projective as a $(B,A)$-module and $M'$ 
is weakly (resp. strongly) $\S$-locally projective as an $(A,B)$-module;
\item $S$ is a $B$-ring with (resp. idempotent) local units, $M$ is a firm 
right $S$-module and $M'$ is a firm left $S$-module;
\item $S$ has left (resp. idempotent) local units on $M$ and right (resp. 
idempotent) local units on $M'$ as a $B$-ring.
\end{enumerate}
In the strong/idempotent case, these conditons are also equivalent 
the following condition: 
\begin{enumerate}
\item[(4)] $M$ is the direct limit of
a split direct system $(P_i)_{i\in I}$, with all the $P_i$ finitely generated and projective as right $A$-modules and such that $M=\varinjlim P_i$ and $M'=\varinjlim P_i^*$.
\end{enumerate}
\end{corollary}

\begin{proof}
$\ul{1.\Rightarrow 2.}$ It follows from \thref{3.6} that $\psi$ is injective,
and then from the final statement in \thref{3.7},
with $M$ replaced by $M'$ and vice versa, that $\Phi$ is injective.
In a similar way, it follows from \thref{3.6}, now with $M$ replaced by $M'$
that $\varphi$ is injective, and then from the final statement in \thref{3.7}
that $\Psi$ is injective. Then the implication $\ul{1.\Rightarrow 2.}$
follows from Theorems \ref{th:3.6} and \ref{th:3.7}.\\
$\ul{2.\Rightarrow 3.}$ and $\ul{3.\Rightarrow 1.}$ follow immediately
from \thref{3.7} (and its analogous version, with the roles of $M$ and
$M'$ interchanged), taking into account the fact that the hypothesis that
$\Phi$ is injective is only needed in the proof of $\ul{1.\Rightarrow 2.}$
in \thref{3.7}.\\

$\ul{4\Rightarrow 1}.$ Since $\varphi$ is injective, we can view
$M'$ as a submodule of $M^*$, and therefore $\R=M'$. For a finitely
generated submodule $N$ of $M$, we can find $i\in I$ such that
$N\subset P_i$. The dual basis of $P_i$ is contained in
$P_i\otimes_A P_i^*\subset M\otimes_A \R$, and is an idempotent dual basis on $N$;
in a similar way, we can find idempotent local bases on finitely generated
submodules of $M'$.\\

$\ul{1 \Rightarrow 4}.$ 
Using \thref{2.17}, we find split direct systems $(P_i)_{i\in I}$ and $(P'_j)_{j\in J}$ such that $M=\varinjlim P_i$ and $M'=\varinjlim P'_j$. 
Moreover, it follows from the proof of \thref{2.17} that 
$P_i=e_iM$, where $e_i\in S$ is a idempotent local dual basis for $P_i$. 
Analogously, $P'_j=M'e'_j$ with $e'_j\in S$ local dual bases for $P'_j$. We will construct from these direct systems a new split direct system that satisfies the desired properties.

Let $K\subset S$ be the set consisting of all the previously considered dual bases $e_i$ and $e'_j$ for respectively $P_i$ and $P'_j$. For two elements $k_1,k_2\in K$ we say that $k_1\le k_2$ if and only if $k_1k_2=k_1$ in $S$. Now put $M_k:=k\cdot M$ and $M'_k:=M'\cdot k$ for all $k\in K$, similar as in the proof of \thref{2.17} (part 3 to 1) we find $M=\varinjlim M_k$ and $M'=\varinjlim M'_k$.

Moreover, $M'_i=(M_i)^*$: taking the restriction and corestriction of the morphism $\varphi:M'\to M^*$, we find a morphism $\varphi_k:M_k'\to M_k^*$, $\varphi_k(m')(m)=\mu(m'\otimes_Bm)$ for all $m'\in M_k'$ and $M\in M$. Let us first prove that this morphism is injective.
Suppose $\mu(m'\otimes_Bm)=0$ for all $m\in M_k$. We know by the construction of $M_k$ that $m'=m'\cdot k$. Also by the construction of $M_k$ and $M'_k$, we can find elements $u_i\in M_k$ and $f_i\in M'_i$ such that $k=u_i\otimes_Af_i$. Furthermore, $m'=m'\cdot k=\mu(m'\ot_B u_i)f_i=0$.

Finally, $\varphi_k$ is also surjective. Take $f\in M_k^*$. Then
$$f(m)=f(k\cdot m)= f(\sum_iu_i\mu(f_i\ot_B m))=\sum_if(u_i)\mu(f_i\ot_B m),$$
for all $m\in M_k$, where we denoted again $k=\sum_iu_i\ot f_i$. We conclude that $f=\varphi_k(f(u_i)f_i)$.
\end{proof}

Let $S$ be a ring with idempotent local units. It follows from
\leref{2.10} that $S=\varinjlim Se=\varinjlim eS=\varinjlim eSe$, where 
the limit is taken over the set of idempotents $e$ of $S$. So $S$ is the direct
limit of a system of rings with left units, of a system of rings with right local units and of
a system of rings with two-sided local units.
In the
case where $S$ is as in \coref{3.8}, we have
$M=\varinjlim P_i$ and $M'=\varinjlim P^*_i$. Since the tensor product commutes with direct limits, we find
$$S=M\otimes_A M'=\varinjlim P_i \otimes_A M' = \varinjlim M\otimes_A P^*_i.$$
Here $ P_i \otimes_A M'$ is a ring with left unit $e_i$, with $e_i$ the finite dual basis for $P_i$, and $e_i$ is a right unit of the ring $M\otimes_A P^*_i$.\\
Recall that $P_i\subseteq P_j$ if $i\leq j$. Take $k\geq i,j$. Then
$P_i\otimes_A P^*_j\subseteq P_k\otimes_A P_k^*$, and it follows that 
$$S=M\otimes_A M'=\varinjlim P_i\otimes_A P_j^*=\varinjlim P_k\otimes_A P^*_k,$$
with $P_k\otimes_A P^*_k$ a ring with unit by \thref{3.4}.

We will now look at the situation where the direct limit is a direct sum.
Recall that a ring $R$ has enough idempotents if there exists a set $\{e_i\}_{i\in I}$ of pairwise orthogonal idempotents, such that every element in $R$ admits a finite sum of these idempotents as a two-sided unit, or, equivalently $R=\bigoplus_{i\in I}e_i R=\bigoplus_{i\in I}Re_i$. 
We call $\{e_i\}_I$ a complete family of idempotents for $R$.

\begin{theorem}\thlabel{3.9}
Let $\ul{M}=(M,M',\mu)$ be a dual pair over the rings $A$ and $B$, and $S=
M\ot_A M'$ the associated elementary $B$-ring. The following statements are equivalent
\begin{enumerate}
\item $S$ is a ring with enough idempotents and $M$ and $M'$ are firm,
respectively as a left and right $S$-module;
\item there exist a family of finitely generated and projective right $A$-modules
$\{P_i~|~i\in I\}$ such that
$$M=\bigoplus_{i\in I} P_i~~{\rm and}~~M'=\bigoplus_{i\in I} P^*_i.$$
\end{enumerate}
\end{theorem}

\begin{proof}
$\ul{(1)\Rightarrow (2)}.$ 
$S$ has a complete family of idempotents $\{e_i~|~i\in I\}$, so it has idempotent local units,
and, by \coref{3.8}, $M$ is strongly $\varphi_\mu(M')$-locally projective and $\varphi_\mu$ is injective. For every $i\in I$, $e_iM$ is now a finitely generated and projective direct summand of $M$ (see \thref{2.17}). We claim 
$M=\bigoplus_{i\in I}e_iM$. Indeed, if $e_im=e_jm'$ with $i\neq j$, then $e_im=e_i^2m=e_ie_jm'=0$.

In a similar way, we find that $M'=\bigoplus_{i\in I}M'e_i$. As in the proof
of \coref{3.8}, we show that $M'e_i\cong (e_iM)^*$.\\

$\ul{(2)\Rightarrow (1)}.$ 
Clearly, $M$ and $M'$ are strongly locally projective, so $S$ is a ring with idempotent local units and $M$ and $M'$ are firm $S$-modules. We are done if we can construct
a complete family of idempotents. Let
$e_i=\sum_ju_i^j\otimes_A f_i^j\in P_i\otimes_A P_i^*$ be the finite dual basis for $P_i$.
$f_i^j(e_k^\ell)=0$ for all $j,\ell$ if $i\neq k$, so $e_ie_k=0$, and the $e_i$
form a complete set of idempotents. We then have
$$S=\bigoplus_{i\in I}e_iP_i\otimes_A M' \quad {\rm and} \quad S=\bigoplus_{i\in I} M\otimes_A P^*_ie_i.$$
\end{proof}

\begin{remark}
For sake of completeness, let us observe the following connection to Morita contexts. If $(M,M',\mu)$ is a dual pair, then we have a Morita context $(A,S,M',M,\mu,S)$, which is strict if and only if $\mu$ is surjective. Conversely, every morita context $(A,S,P,Q,f,g)$ between a ring $A$ and a $B$-ring $S$ gives rise to a dual pair $({_BQ_A},{_AP_B},\bar{f})$, where $\bar{f}$ is induced by $f$. The original context is naturaly homomorphic to the context associated to this dual pair. If the map $g$ from the original context is surjective, then both contexts are isomorphic. 
If $S$ has (idempotent) local units, then by \coref{3.8}, $P$ and $Q$ are weakly (strongly) locally projective $A$-modules, if $P$ and $Q$ are firm $S$-modules. By lemma \leref{3.3}, this last condition will be satisfied if $f$ is surjective.
\end{remark}

\section{Local comultiplications}\selabel{4}
Let $\ul{M}=(M,M',\mu)$ be a dual pair over the rings $A$ and $B$,
and consider the $A$-bimodule $\Cc=M'\otimes_BM$. Write
$\varepsilon_\Cc=\mu:\ \Cc\to A$. In order to make $\Cc$ into an $A$-coring,
we need a comultiplication. 
Dualising the construction of local units for the $B$-ring $S=M\ot_AM'$, we can 
start from an element $e\in S^B$ and define 
$\Delta_e(m'\otimes_Bm)=m'\otimes_Be\otimes_Bm\in \C\otimes_A\C$. This map 
is clearly coassociative, however the counit property 
$m_i\otimes_Be\cdot m'_i=c_{(1)}\varepsilon(c_{(2)})=c=\varepsilon(c_{(1)})c_{(2)}=m_i\cdot e\otimes_Bm'_i$ holds only for 
elements $c=m'_i\otimes_Bm_i$ such that $e$ is a local dual basis (local unit) for 
the $m_i$ and $m'_i$. This leads us to \deref{4.1}.

\begin{definition}\delabel{4.1}
Let $A$ be a ring with unit, $\Cc$ an $A$-bicomodule and $\varepsilon_\Cc:\ \C\to A$ an $A$-bimodule
map. A right
\emph{$\varepsilon$-comultiplication} on $\Dd\subset \Cc$ is a 
coassociative $A$-bimodule map
$\Delta_\D:\ \C\to \C\otimes_A\C$ such that
$c=c_{(1)}\varepsilon_\Cc(c_{(2)})=(\Cc\otimes_A\varepsilon_\C)\circ\Delta_\D(c)$ 
for all $c\in\D$. A right
$\varepsilon$-comultiplication $\Delta_\D$ is called \emph{idempotent} if $(\Cc\otimes_A\varepsilon)\circ\Delta_\Dd$ is an idempotent in ${_A\End_A}(\Cc)$.

If there exists a right $\varepsilon$-comultiplication on $\C$, we say that $\Cc$ is an $A$-coring with right comultiplication.

We call $\Cc$ an $A$-coring with \emph{weak} (resp. strong) \emph{right local 
$\varepsilon$-comultiplications} 
if for every finitely generated right $A$-submodule $\Dd$ of $\Cc$, we can find a right $\varepsilon$-comultiplication 
(resp. an idempotent right $\varepsilon$-comultiplication) $\Delta_\Dd$ on $\Dd$.

In a similar way, we define corings with weak and strong left local
$\varepsilon$-comultiplications.

We say that $\Cc$ is an $A$-coring with two-sided
weak (resp. strong) right local $\varepsilon$-comultiplic\-ations if
for every finitely generated right $A$-submodule $\D\subset\C$ there 
exists an $A$-bilinear coassociative map $\Delta_\D$ which is at the same time a 
right and a left $\varepsilon$-comultiplication
(resp. a strong right and a left $\varepsilon$-comultiplication). 

A right $A$-module $M$ is called a \emph{weak local right $\C$-comodule} if for every finitely generated $A$-submodule $N\subset M$, there exists a comultiplication
$\Delta_N:\ \C\to\C\otimes_A\C$ on $\Cc$ and a right $A$-linear map
$\delta_N:\ M\to M\otimes_A\C, \delta_N(m)=m_{[0]}\ot_A m_{[1]}$ 
such that $n_{[0]}\varepsilon(n_{[1]})=n$ for all $n\in N$ and
$\delta_N(m_{[0]})\otimes_Am_{[1]}=m_{[0]}\otimes_A\Delta_N(m_{[1]})$ for
every $m\in M$. We call $M$ a \emph{strong local right $\C$-comodule} if $M$ is a 
weak local right $\C$-comodule and, in addition, $(\Cc\otimes_A\varepsilon)\circ\delta_N$ is an idempotent in ${\End_A}(M)$, for every finitely generated $N\subset M$.

We say that two comultiplications $\Delta$ and $\Delta'$ on $\C$  
\emph{coassociate} if
$(\C\otimes_A \Delta)\circ\Delta'=(\Delta'\otimes_A \C)\circ\Delta$ 
and
$(\C\otimes_A \Delta')\circ\Delta=(\Delta\otimes_A \C)\circ\Delta'$.
\end{definition}

Observe that the fact that $\C$ is a local comodule over $\C$ does not imply 
that $\C$ is an $A$-coring with local comultiplications, since it is possible that $\delta_\D\neq\Delta_\D$.

\begin{theorem}\thlabel{4.2}
If $\Cc$ is an $A$-coring with (left, right) weak (resp. strong) local comultiplications then ${_A\End_A}(\Cc)$ has (left, right) (resp. idempotent) local units on $\Cc$.
\end{theorem}

\begin{proof}
The local units are of the type $(\Cc\otimes\varepsilon)\circ\Delta_\Dd$, where $\Delta_\Dd$ is a right local comultiplication on a finitely generated submodule $\Dd\subset\Cc$.
\end{proof}

\begin{theorem}\thlabel{4.3}
Let $\C$ be an $A$-bimodule and $\varepsilon_\Cc:\ \C\to A$ an $A$-bimodule 
map. Then the following statements are equivalent.
\begin{enumerate}
\item For every $c\in\C$, there exists a right 
$\varepsilon$-comultiplication $\Delta_c$ on $\{c\}$, 
such that $\Delta_c$ 
and $\Delta_{c'}$ coassociate,
for all $c,c'\in \C$;
\item there exists a right $\varepsilon$-comultiplication on every 
$A$-subbimodule of 
 $\C$ generated by a single element such that two 
such comultiplications 
coassociate;
\item $\C$ is an $A$-coring with right weak local comultiplications.
\end{enumerate}
\end{theorem}

\begin{proof}
$\underline{1.\Rightarrow 2.}$ Take $c\in\C$ 
and $\Delta_c$ as in part 1. We only have to prove that the 
counit property holds for every element of the form $acb$, with $a,b\in 
A$. But $\Delta_c(acb)=ac_{(1)}\otimes_Ac_{(2)}b$, so 
$ac_{(1)}\varepsilon(c_{(2)}b)=ac_{(1)}\varepsilon(c_{(2)})b=acb$.\\

$\underline{2.\Rightarrow 3.}$ Let $\D$ be the  
$A$-subbimodule of $\C$ generated by the elements $\{c^1,\ldots,c^k\}$.
We proceed by 
induction on the number of generators. Let $\Delta$ be the 
$\varepsilon$-comultiplication on the $A$-subbimodule generated by 
$c^1$, and $\Delta'$ the $\varepsilon$-comultiplication on the $k-1$ elements 
$c^i-c^i_{(1)}\varepsilon(c^i_{(2)})$, $i=2,\ldots,k$. By assumption, 
these comultiplications can be chosen in such a way that they coassociate.
Now  
$\Delta''=\Delta+\Delta'-\Delta'\circ(\C\otimes_A\varepsilon)\circ\Delta$
is a comultiplication on $\C$. It is obvious that $\Delta''$ is an
$A$-bilinear map; let us sketch the proof of the coassociativity. For bimodule
maps $\Delta$ and $\varepsilon$,  we always have, without any counit
property assumption, that
$$\Delta\circ(\varepsilon\otimes_A\C)=
(\varepsilon\otimes_A\C\otimes_A\C)\circ(\Cc\otimes_A\Delta).$$ 
Using this property and the (mixed) coassociativity of $\Delta$ and $\Delta'$, 
we can show that $\Delta''$ is coassociative. We restrict ourselves to 
proving the following identity, leaving all other details to the reader.
\begin{eqnarray*}
&&\hspace*{-2cm}(\Cc \ot_A \Cc \ot_A \Cc \ot_A\Delta)\circ(\Cc \ot_A \Cc \ot_A\Delta')
\circ(\Cc \ot_A\Delta)\circ\Delta'\\
&=&(\Cc \ot_A \Cc \ot_A \Delta'\ot_A \Cc )\circ(\Cc \ot_A \Cc \ot_A\Delta)
\circ(\Cc \ot_A\Delta)\circ\Delta'\\
&=&(\Cc \ot_A \Cc \ot_A \Delta'\ot_A \Cc )\circ(\Cc \ot_A\Delta\ot_A \Cc )
\circ(\Cc \ot_A\Delta)\circ\Delta'\\
&=&(\Cc \ot_A \Delta\ot_A \Cc \ot_A \Cc )\circ(\Cc \ot_A\Delta'\ot_A \Cc )
\circ(\Cc \ot_A\Delta)\circ\Delta'\\
&=&(\Cc \ot_A \Delta\ot_A \Cc \ot_A \Cc )\circ(\Cc \ot_A\Delta'\ot_A \Cc )
\circ(\Delta'\ot_A \Cc )\circ\Delta\\
&=&(\Cc \ot_A \Delta\ot_A \Cc \ot_A \Cc )\circ(\Delta'\ot_A \Cc \ot_A \Cc )
\circ(\Delta'\ot_A \Cc )\circ\Delta\\
&=&(\Delta'\ot_A \Cc \ot_A \Cc \ot_A \Cc )\circ(\Delta\ot_A \Cc \ot_A \Cc )
\circ(\Delta'\ot_A \Cc )\circ\Delta.
\end{eqnarray*}
We now check that the counit property holds for all elements in $\D$. 
We know that $(\Cc\otimes_A\varepsilon)\circ\Delta(c^1)=c^1$, and since
$\Delta''(c^1)=\Delta(c^1)+\Delta'(c^1)-\Delta'(c^1)$, this implies that $\Delta''$ is 
also an $\varepsilon$-comultiplication on $c^1$. Furthermore, for every 
$i=2,\ldots, k$ we have
\begin{eqnarray*}
&&\hspace{-2cm}
(I\otimes\varepsilon)\circ\Delta''(c^i)=
(I\otimes\varepsilon)\circ\Delta(c^i)+
(I\otimes\varepsilon)\circ\Delta'\circ(c^i-c^i_{(1)}\varepsilon(c^i_{(2)}))\\
&=&c^i_{(1)}\varepsilon(c^i_{(2)})+c^i-c^i_{(1)}\varepsilon(c^i_{(2)})=c^i,
\end{eqnarray*}
and we conclude that $c=c_{(1'')}\varepsilon(c_{(2'')})$ for all 
$c\in\D$.\\

$\underline{3.\Rightarrow 1.}$ We know now 
 that there exists an $\varepsilon$-comultiplication for every $c\in\C$. Given two 
elements, we know that there exists a common 
$\varepsilon$-comultiplication, so the coassociativity is 
automatically satisfied.
\end{proof}

\begin{theorem}\thlabel{4.4}
If there exists a left $\varepsilon$-comultiplication $\Delta$ on $c\in 
\C$ and a right $\varepsilon$-comultiplication $\Delta'$ on $c'\in\C$ 
such that $\Delta$ and $\Delta'$ coassociate, 
then there exists a coassociative $A$-bimodule map $\Delta''$ which is 
left $\varepsilon$-comultiplication on $c$ and a right 
$\varepsilon$-comultiplication on $c'$ 
\end{theorem}

\begin{proof}
Put $\Delta''=\Delta+\Delta'-
(\Cc\otimes_A\varepsilon\otimes_A\Cc)\circ(\Delta'\otimes_A\Cc)\circ\Delta$. The 
coassociativity of $\Delta''$ can be proven along the same lines as in the 
previous theorem. Using the bilinearity of $\varepsilon$ and $\Delta'$ we 
find the following identity
\begin{eqnarray*}
&&\hspace*{-2cm}(\varepsilon\otimes_A \Cc )\circ(\Cc \otimes_A \varepsilon\otimes_A\Cc )\circ
(\Delta'\otimes_A\Cc )\circ\Delta\\
&=&(\varepsilon\otimes_A \Cc )\circ(\varepsilon\otimes_A \Cc \otimes_A \Cc )\circ
(\Delta'\otimes_A\Cc )\circ\Delta\\
&=&(\varepsilon\otimes_A \Cc )\circ\Delta'\circ
(\varepsilon\otimes_A \Cc )\circ\Delta.
\end{eqnarray*}
We now apply this to prove the left counit property on $c$.
\begin{eqnarray*}
&&\hspace*{-2cm}
(\varepsilon\otimes_A\Cc )\circ\Delta''(c)=
(\varepsilon\otimes_A\Cc )\circ\Delta(c)+
(\varepsilon\otimes_A\Cc )\circ\Delta'(c)\\
& &-(\varepsilon\otimes_A\Cc )\circ
(\Cc \otimes_A\varepsilon\otimes_A\Cc )\circ(\Delta'\otimes_A\Cc )\circ\Delta(c)\\
&=&c+(\varepsilon\otimes_A\Cc )\circ\Delta'(c)-
(\varepsilon\otimes_A \Cc )\circ\Delta'\circ
(\varepsilon\otimes_A \Cc )\circ\Delta(c)\\
&=&c+(\varepsilon\otimes_A\Cc )\circ\Delta'(c)-
(\varepsilon\otimes_A \Cc )\circ\Delta'(c)=c
\end{eqnarray*}
Using the coassociativity, one proves in a similar way that $\Delta''$ is a
right $\varepsilon$-comulti-plication on $c'$.
\end{proof}

\begin{corollary}\colabel{4.5}
If $\Cc$ is an $A$-coring with left and right weak local comultiplications, then $\Cc$ has also two-sided weak local comultiplications.
\end{corollary}

\begin{theorem}\thlabel{4.6}
Let $\Cc$ be an $A$-bimodule and $\varepsilon :\ \Cc\to A$ an $A$-bimodule map. The following statements are equivalent
\begin{enumerate}
\item
$\Cc$ is an $A$-coring with right strong local comultiplications; 
\item every finitely generated right $A$-submodule $\Dd\subset\Cc$ is contained in a $A$-subbimodule $\Ee\subset\Cc$, which is an $A$-coring with right comultiplication and a direct summand of $\Cc$,
\item 
there exists a split direct system $(\Cc^i)_{i\in I}$, where $\Cc^i$ is an $A$-coring with a right comultiplication, such that $\Cc=\varinjlim \Cc^i$.
\end{enumerate}
\end{theorem}

\begin{proof}
$\ul{1\Rightarrow 2.}$
Suppose $\Cc$ has right strong local comultiplications and let $\Dd$ be a finitely generated right $A$-submodule of $\Cc$. Then we know $\psi=(\Cc\otimes_A\varepsilon) \circ \Delta_\Dd$ is an idempotent in ${_A\End_A(\Cc)}$ and thus a projection. Denote $\Ee:=\im \psi$ and define a new comultiplication
$$\Delta_\Ee=(\Cc\ot_A\varepsilon\ot_A\Cc)\circ(\Delta_\Dd\ot_A\Delta_\Dd):\Cc\to\Cc\ot_A\Cc$$
It is an easy computation to check that this map is coassociative. Furthermore, for $e=c_{(1)}\varepsilon(c_{(2)})\in\Ee$, with $c\in\Cc$, we find 
$\Delta_\Ee(e)=c_{(1)}\varepsilon(c_{(2)})\ot_Ac_{(3)}\varepsilon(c_{(4)})\in\Ee\ot_A\Ee.$
Moreover, 
\[
\begin{array}{rcl}
(\Cc\ot_A\varepsilon)\circ\Delta_\Ee(e)&=&
(\Cc\ot_A\varepsilon)\circ(\Cc\ot_A\varepsilon\ot_A\Cc)\circ(\Delta_\Dd\ot_A\Delta_\Dd)
\circ(\Cc\ot_A\varepsilon)\circ\Delta_\Dd (c)\\
&=&(\Cc\ot_A\varepsilon)\circ\Delta_\Dd \circ(\Cc\ot_A\varepsilon)\circ\Delta_\Dd 
\circ(\Cc\ot_A\varepsilon)\circ\Delta_\Dd (c)\\
&=&(\Cc\ot_A\varepsilon)\circ\Delta_\Dd (c)=e
\end{array}
\]
We can conclude that $\Ee$ is an $A$-coring with right comultiplication. Finally $\Dd\subset\Ee$, since $c=\psi(c)$ for all $c\in\Dd$.

$\ul{2\Rightarrow 3.}$ 
Denote by $I$ the set consisting of all $A$-corings
with right comultiplication that are direct summands of $\Cc$. This set is
partially ordered: for $\Ee, \Ee'\in I$ we define 
$\Ee\le \Ee'$ if $\Delta_{\Ee'}$ is a right $\varepsilon$-comultiplication on $\Ee$. In this situation, the projection $\psi_\Ee=(\Cc\ot_A\varepsilon)\circ\Delta_\Ee$ factors trough $\psi_{\Ee'}$, and the rest follows easily.

$\ul{3\Rightarrow 1.}$
For every finitely generated right $A$-submodule $\Dd\subset\Cc$, we can find a $\Cc^i$ containing $\Dd$. Since $\Cc^i$ is a direct summand of $\Cc$, we can extend the comultiplication on $\Cc^i$ to the whole of $\Cc$ by making it zero on the complement. This is a right $\varepsilon$-comultiplication on $\Dd$ and this finishes the proof.
\end{proof}

\begin{theorem}\thlabel{4.7}
Let $\Cc$ be an $A$-bimodule and $\varepsilon :\ \Cc\to A$ an $A$-bimodule map. The following statements are equivalent
\begin{enumerate}
\item
$\Cc$ is an $A$-coring with two-sided strong local comultiplications; 
\item every finitely generated $A$-bimodule $\Dd\subset\Cc$ is contained in a $A$-bisubmodule $\Ee\subset\Cc$, such that $\Ee$ is an $A$-coring and a direct summand of $\Cc$;
\item 
there exists a split direct system $(\Cc^i)_{i\in I}$ of $A$-corings, such that $\Cc=\varinjlim \Cc^i$.
\end{enumerate}
\end{theorem}

\begin{proof}
$\ul{1\Rightarrow 2.}$ Let $\Delta$ be the left and right $\varepsilon$-comultiplication of a finitely generated $A$-bisubmodule $\Dd\subset\Cc$. 
Since $\Cc$ has two-sided strong local comultiplications,
$\alpha=(\Cc\otimes_A\varepsilon)\circ\Delta$ and $\beta=(\varepsilon\otimes_A\Cc)\circ\Delta$ are idempotents. $\alpha$ and
$\beta$ commute, since $(\alpha\circ\beta)(c)=\varepsilon(c_{(1)})c_{(2)}
\varepsilon(c_{(3)})=(\beta\circ\alpha)(c)$, for all $c\in \Cc$.
Therefore $\gamma=\alpha\circ\beta=\beta\circ\alpha$ is also an
idempotent, and a projection $\Cc\to \Cc$. Let $\Ee$ be the image of
$\gamma$, then the restrictions of $\alpha$ and $\beta$ to $\Ee$
are the identity, since $\alpha\circ\gamma=\gamma$ and $\beta\circ\gamma=\gamma$. We obtain that $\Ee$ is an $A$-coring with comultiplication $\Delta$.
The proof of the other implications is similar to the one of the corresponding
implications in \thref{4.6}.
\end{proof}

\begin{blanco}{\bf Local projectivity versus local comultiplications.}
Take a dual pair 
$({_BM_A},{_AM'_B},\mu)$, and assume
that the associated elementary algebra $S=M\ot_A M'$ is a $B$-ring with local units 
on $M$ and $M'$ (which is equivalent to $M$ being 
$M'$-locally projective and $M'$ being $M$-locally projective, in such a way that 
the dual bases can be chosen in $S^B$). Then 
$\Cc=M'\otimes_BM$ is an $A$-coring with weak local comultiplications and $M$ and $M'$ are local $\Cc$-comodules.
If the local units can be chosen to be idempotent, then $\C$ is a coring 
with strong local comultiplications. Observe that we do not have the
converse implication: the fact that $M$ and $M'$ are local 
comodules, does not necessarily imply that $M$ and $M'$ have local dual bases.

As a special case of this construction, we recover the construction of
 a comatrix coring (see \cite{EGT}): let $M\in {_B\M_A}$ be finitely
 generated and projective as a right $A$-module,
and consider the dual pair $(M,M^*,\mu)$ where $\mu(f\otimes_Bm)=f(m)$ for all 
$m\in M$ and $f\in M^*$. We 
can find an idempotent dual basis $e=u_i\otimes_A f_i\in S$. Since this 
is a dual basis 
for every element in $M$, it is a dual basis for every $B$-submodule of 
$M$. It follows from \coref{3.8} that $e\in S^B$. 
This means 
we can construct a local $\varepsilon$-comultiplication on every finitely 
generated 
submodule of $\C$, and since $M$ itself is finitely generated, we have a 
usual comultiplication and $\Cc$ is an $A$-coring. The comultiplication and 
counit are explicitly defined by
\begin{equation}\eqlabel{4.7.1}
\Delta_\C(f\otimes_B m)=\sum_if\otimes_B u_i\otimes_A f_i\otimes_B m
, \qquad \varepsilon_\C(f\otimes_B m)=f(m),
\end{equation}
for all $f\otimes_B m\in \C$. The conditions in \thref{3.4} are then
equivalent to
\begin{itemize}
\item[(7)] There exists $e\in (M\ot_A M')^B$ such that \equref{4.7.1}
defines an $A$-coring structure on $\Cc=M'\ot_B M$, and $M$ and $M'$
are respectively a right and left $\Cc$-comodule, with coactions
$$\rho^r(m)=e\ot_B m~~{\rm and}~~\rho^l(m')=m'\ot_B e.$$
\end{itemize}
\end{blanco}

\section{Local counits}\selabel{5}
\begin{definition}\delabel{5.1}
Let $A$ be a ring with unit,
$\Delta_\Cc$  a coassociative comultiplication on an $A$-bimodule $\Cc$,
and $M$ a right $A$-module with a coassociative right $A$-coaction
$\delta_M:\ M\to M\ot_A\Cc$. A right \emph{counit on $M$} is an $A$-bimodule
map $\varepsilon_M:\ \Cc\to A$ such that
$$(M\ot_A \varepsilon_M)\circ\delta_M=M.$$
$\varepsilon_M$ is called idempotent if $(\varepsilon_M\ot_A\varepsilon_M)\circ\Delta_\Cc
=\varepsilon_M$.\\
We say that $M$ has (idempotent) right  \emph{local counits} if there exists an (idempotent)
right counit $\varepsilon_N$ on every finitely generated
right $\Cc$-subcomodule $N\subset M$.\\
If there exist right (idempotent) local counits on $\Cc$, then we call
$\Cc$ an $A$-coring with right (idempotent) local counits.\\
Left and two-sided (idempotent) local counits can be introduced in a similar
way.
\end{definition}

The terminology ``idempotent" counit is justified by the following Lemma.

\begin{lemma}\lelabel{5.2}
Let $\varepsilon_M$ be a  right counit on a right $\C$-comodule $M$.
The following are equivalent.
\begin{enumerate}
\item $\varepsilon_M$ is an idempotent counit;
\item $\varepsilon_M:\ \Cc\to A$ is comultiplicative (the comultiplication
of $A$ is the canonical isomorphism $A\cong A\ot_AA$);
\item $\varepsilon_M$ is an idempotent in $\*C$.
\end{enumerate}
\end{lemma}

\begin{proof}
The equivalence of (1) and (2) is obvious. In $\*C$, we easily compute that $(\varepsilon_M*\varepsilon_M)(c)=
\varepsilon_M(c_{(1)}\varepsilon_M(c_{(2)}))=\varepsilon_M(c_{(1)})
\varepsilon_M(c_{(2)})$, and we deduce immediately the equivalence of (1) and (3).
\end{proof}

\begin{theorem}\thlabel{5.3}
There exists (idempotent) right local counits on a right 
$\C$-comodule $M$ if and only if the $A$-ring $\*C$ has (idempotent) right local units on $M$.
In particular, $\C$ is an $A$-coring with (idempotent) right local counits if
and only if the $A$-ring $\*C$ has (idempotent) right local units on $\C$.
\end{theorem}

\begin{proof}
The right counits on $M$ are precisely the units in $\*C$ on $M$.
\end{proof}

The results of \seref{2.1} can be restated in terms of local counits.
In particular, the existence of a counit on every $c\in\Cc$ implies 
the existence of local counits on $\C$, and if $\Cc$ has (left, right or two-sided) idempotent local counits, then we can write $\Cc=\varinjlim \Cc^i$, where $(\Cc^i)_{i\in I}$ is a split direct system of $A$-corings with a (left, right or two-sided) counit.

We call $M\in\Mm^\Cc$ \emph{cofirm} if the right coaction $\delta_M$ on $M$ induces
an isomorphism $M\cong M\otimes^\Cc\Cc$ in $\Mm^\Cc$. 

\begin{theorem}\thlabel{5.4}
If $\Cc$ is an $A$-coring with right local counits, a right $\Cc$-comodule $M$
is cofirm if and only if $\Cc$ has right local counits on $M$.
\end{theorem}

\begin{proof}
First assume that $\Cc$ has right local counits on $M$. Suppose that $\delta_M(m)=0$.
Take a local counit $\varepsilon_m$ on $m$; then $m=(M\ot\varepsilon_m)(\delta_M(m))=0$,
and it follows that $\delta_M$ is injective. Take 
$\sum_i m_i\ot_A c_i\in M\ot^\Cc \Cc$, and a right local counit $\varepsilon$
on the right $A$-submodule of $\Cc$ generated by $\{c_1,\cdots,c_n\}$. Then we
compute that
$$\delta_M(\sum_i m_i\varepsilon(c_i))=
\sum_i m_{i[0]}\ot_A m_{i[1]}\varepsilon(c_i)=
\sum_i m_{i}\ot_A c_{i(1)}\varepsilon(c_{i(2)})=\sum_i m_i\ot_A c_i,$$
so it follows that $\delta_M:\ M\to M\ot^\Cc\Cc$ is surjective.\\
Conversely, if $\Cc$ is an $A$-coring with right local counits, $\Cc$ has right local counits on $M\otimes^\Cc\Cc$. If $M$ is cofirm, then $M\otimes^\Cc\Cc\cong M$, hence
$\Cc$ also has right local counits on $M$.
\end{proof}

\section{The Dorroh coring}\selabel{6}
Let $B$ be a ring, and $R$ a $B$-ring without unit. Recall the construction
of the so-called Dorroh overring: $S=R\times B$ is a $B$-bimodule, with the following
left and right $B$-action: $b'(r,b)b''=(b'rb'',b'bb'')$. The map
$\eta:\ B\to R\times B$, $\eta(b)=(0,b)$ is a $B$-bimodule map. $S$ is a
$B$-ring with unit map $\eta$ and multiplication
$(r,b)\cdot(r',b')=(rr'+rb'+r'b,bb')$. We refer to e.g. to \cite[section 1.5]{W2}
in the case where $B=\ZZ$.\\
We will now present a dual version of this construction, for corings. 

Recall that 
an $A$-subbimodule of a counital $A$-coring $\Cc$ is called a (two-sided) 
coideal if $\Delta(\D)\subseteq 
\C\otimes_A\D + \D\otimes_A\C$ and $\varepsilon_\C(\D)=0$.

\begin{theorem}\thlabel{6.1}
Let $\Cc$ be an $A$-coring, not necessarily with a counit. There exists a counital 
$A$-coring $\hat{\C}$ with the following properties:
\begin{enumerate}
\item $\Cc$ is isomorphic to a coideal of $\hat{\C}$;
\item there exists a surjective $A$-coring morphism $\pi:\ \hat{\C}\to \C$;
\item there exists an injective $A$-coring morphism $\iota:\ A\to \hat{\C}$;
\item the category of (not necessarily counital) comodules over $\Cc$
is isomorphic to the category of counital comodules over $\hat{\Cc}$.
\end{enumerate}
We call $\hat{\Cc}$ the Dorroh coring associated to $\Cc$.
\end{theorem}

\begin{proof}
$\hat{\C}=\Cc\times A$ is an $A$-bimodule with left and right $A$-action given by 
$b\cdot(c,a)b'=(bcb',bab')$, for all $c\in\Cc$ and $a,b,b'\in A$. 
The map $\varepsilon:\ \hat{\C}\to A$, 
$\varepsilon(c,a)=a$ is an $A$-bimodule map. It is easy to verify that
$\hat{\Delta}:\ \hat{\C}\to \hat{\C}\ot_A\hat{\C}$, given by the formula
$$\hat{\Delta}(c,a)=(c_{(1)},0)\otimes_A(c_{(2)},0) 
+ (0,1)\otimes_A(c,a) + (c,a)\otimes_A(0,1)-(0,a)\otimes_A(0,1),$$
is a coassociative $A$-bimodule map, and that $(\hat{\Cc},\hat{\Delta},\varepsilon)$
is a counital coring. The map $\pi$ is the canonical surjection $\hat{\Cc}\to \Cc$,
and $\iota$ is the canonical injection $A\to \hat{\Cc}$. All further verifications are
straightforward.\\
Let $(M,\delta_M)$ be a right $\Cc$-comodule. A $\hat{\Cc}$-coaction $\hat{\delta}_M$ on $M$ is
defined as follows:
$$\hat{\delta}_M:\ M\to M\otimes_A\hat{\C},~~ 
\hat{\delta}_M(m)=m_{[0]}\otimes_A(m_{[1]},0)+m\otimes_A(0,1),$$
for all $m\in M$. It is straightforward to show that $(M, \hat{\delta}_M)$
is a right $\hat{\Cc}$-comodule. Conversely, if $(M, \hat{\delta}_M)$
is a right $\hat{\Cc}$-comodule, then $(M, \delta_M=(M\otimes_A\pi)\circ\hat{\delta}_M)$
is a right $\Cc$-comodule. Both constructions are functorial, and define
an isomorphism between the categories of right $\Cc$-comodules and
unital right $\hat{\Cc}$-comodules.
\end{proof}

Remark that $\Cc$ is not a subcoring of $\hat{\Cc}$, since 
$\hat{\Delta}(\C)$ is not included in $\C\otimes_A\C$. This has to be
compared with the fact that the quotient map from the Dorroh overring to the original
ring is not a ring morphism.

Finally remark that if $C$ is a coalgebra over a commutative ring $k$, the Dorroh coring $\hat{C}$ is also a coalgebra. Moreover, in this situation, if $C$ is cocommutative, than $\hat{C}$ is also cocommutative.

\section{Rings with local multiplication}\selabel{7}
Let $B$ be a ring with unit, $R$ a $B$-bimodule, and $\eta:\ B\to R$
a $B$-bimodule map. A right \emph{$\eta$-multiplication} on $T\subset R$
is an associative $B$-bimodule map $\mu_T:\ R\otimes_BR\to R$ such that the 
following diagram commutes on the image of $i$:
\[
\xymatrix{
0\ar[r] & T \ar[r]^i & R \ar[r]^{\cong~~~} & R\otimes_BB 
\ar[d]^{I_R\otimes_B\eta}\\
& & &R\otimes_BR \ar[ul]^{\mu_T} }
\]
If $R$ has a right $\eta$-multiplication on itself, then we call $R$ an $A$-ring with right multiplication.

We call $R$ a ring with \emph{weak right local $\eta$-multiplications} if
there exists a right $\eta$-multiplication on every finitely generated
$B$-subbimodule of $R$.

Left and two-sided $\eta$-multiplications are defined in a similar way.

We say that $R$ has \emph{strong right local $\eta$-multiplications} if for every 
finitely generated $B$-subbimodule $T\subset R$, there is a $B$-subbimodule $S\subset R$ containing $T$ on which there exists  a right $\eta$-multiplication $\mu_S$, such that $\mu(R\otimes_B R)\subseteq S$.

Let $M$ be a right $B$-module. 
We say that $R$ has weak right local multiplications on  $M$ if for every finitely generated $B$-submodule $N\subset M$, there exist
a right $B$-linear map $\nu:\ M\otimes_BR\to M$ and a $B$-bimodule map
$\mu:\ R\otimes_BR\to R$ satisfying the usual associativity conditions, and such that $\nu(n\otimes_B\eta(1_B))=n$ for all $n\in N$.

$R$ has strong right local multiplications on $M$ if for every finitely generated $B$-submodule $N\subset M$, we can find a $B$-submodule $N'\subset M$, containing $N$, together with a right $B$-linear map $\nu:\ M\otimes_BR\to M$ and a $B$-bimodule map
$\mu:\ R\otimes_BR\to R$ satisfying the usual associativity conditions, and such that
$\nu(n\otimes_B\eta(1_B))=n$ for all $n\in N'$, and, in addition,
 $\nu(M\otimes_B R)\subset N$.
 
 We say that
 two multiplications $\mu$ and $\mu'$ on $R$  
\emph{associate} if $\mu\circ(R\otimes_B\mu')=\mu'\circ(\mu\otimes_B R)$ 
and $\mu'\circ(R\otimes_B\mu)=\mu\circ(\mu'\otimes_B R)$.\\

Remark that $\eta$ is completely determined by 
$\eta(1)=e\in R^B$, and an associative map $\mu_T:\ R\ot_B R\to R$ is a right $\eta$-multiplication on 
$T\subset R$ if and only if $\mu(r\otimes_B e)= r$ for all $r\in T$.\\
For $b\in B$ and $r\in R$, we will use the notation $\mu_T(r\otimes 
\eta(b))=\mu(r\otimes b)$.

\begin{theorem}\thlabel{7.1}
Let $R$ be a $B$-bimodule and $\eta:\ B\to R$ a $B$-bimodule map. The 
following statements are equivalent.
\begin{enumerate}
\item For every $r\in R$, there exists a right $\eta$-multiplication 
$\mu_r$ on $\{r\}$, such that $\mu_r$ and 
$\mu_{r'}$ associate, for all $r,r'\in R$;
\item there exists a right $\eta$-multiplication on every $B$-subbimodule 
of $R$ generated by a single element such that the $\eta$-multiplications on 
two such subbimodules associate;
\item $R$ is a ring with right weak local $\eta$-multiplications.
\end{enumerate}
\end{theorem}

\begin{proof} 
$\underline{1.\Rightarrow 2.}$ Let $r$ be the generator of a 
$B$-subbimodule $M$ of $R$. We know that there exists a right 
$\eta$-multipication $\mu_r$ on $\{r\}$. If we denote $\mu_r(s,t)=s\cdot 
t$, for all $s,t\in R$, we just verify that $arb\cdot e= ar \cdot 
(be)=ar\cdot (eb)= a(r\cdot e)b=arb$, for all $a,b \in B$, so $\mu$ is a 
right $\eta$-multiplication on $M$.

$\underline{2.\Rightarrow 3.}$ Let $T$ be a finitely generated
$B$-subbimodule of $R$ and $\{t_1,\ldots,t_k\}$ a set of generators for
$T$. We proceed by induction on $k$. Let $\mu$ be a multiplication on the
$B$-bimodule generated by $\{t_2,\ldots, t_n\}$ and denote
$\mu(r,s)=r\cdot s$ for all $r,s\in R$. Let $\mu'$ be a multiplication on
the $B$-bimodule generated by $t_1-t_1\cdot e$ and choose $\mu$ and $\mu'$
in such a way that they associate. Denote $\mu'(r,s)=r*s$. Now define a new multiplication on
$R$ by $r\circ s=r\cdot s +r* s -r\cdot 1 * s$. We leave it to the reader 
to verify that this is an associative right $\eta$-multiplication on $T$.

$\underline{3.\Rightarrow 1.}$ is trivial.
\end{proof}

\begin{theorem}\thlabel{7.2}
If there exists a right $\eta$-multiplication on $t\in R$ and a left 
$\eta$-multiplication on $s\in R$, such that $\mu$ and $\mu'$ associate,  
then there exists a multiplication $\mu''$ on $R$ which is a 
right $\eta$-multiplication on $t$ and a left $\eta$-multiplication on $s$
\end{theorem}

\begin{proof} 
Write $\mu(r,r')=r\cdot r'$ and $\mu'(r,r')=r*r'$, and define $\mu''(r,r')=r\cdot r'+r*r'-r\cdot 
1*r'$.
\end{proof}

\begin{corollary}\colabel{7.3}
$R$ is a ring with left and right local weak $\eta$-multiplications if
and only if $R$ is 
a ring with two-sided weak local $\eta$-multiplications.
\end{corollary}

\begin{theorem}\thlabel{7.4}
Let $R$ be a $B$-bimodule and $\eta:\ B\to R$ a $B$-bimodule map.
The following statements are equivalent.
\begin{enumerate}
\item $R$ is a ring with strong right local multiplications;
\item every finitely generated $B$-subbimodule of $R$
is contained in a direct summand $S$ of $R$, which is a $B$-ring with right multiplication;
\item $R$ can be written as the direct limit of a split direct system $(S_i)_{i\in I}$ of $B$-rings with right multiplication.
\end{enumerate}
\end{theorem}

\begin{proof}
$\ul{1\Rightarrow 2}.$ As in the definition, consider for a $B$-subbimodule
$T$ of $R$ a subbimodule $S$ of $R$ containing $T$, on which there exists
a right $\eta$-multiplication with $\mu(R\otimes_BR)\subset S$.
Then $\mu(S\otimes_BS)\subset S$. We also find that $e=\mu(\eta(1)\otimes_B\eta(1))\in S$. Now we define a new $B$-bilinear map $\eta':\ B\to S$ by $\eta'(b)=eb=be$. It is easy to see that $S$ is a $B$-ring with right multiplication and unit $e$. The map $\psi_S:\ R\to S$, $\psi_S(r)=\mu(r\otimes e)$ is a projection onto $S$, 
and $S$ is a direct summand of $R$.

$\ul{2\Rightarrow 3}$. Take the split direct system of all $B$-rings with right multiplication $S$ as constructed in the previous part and define $S\le S'$ if $\mu_{S'}$ is a right $\eta$-multiplication on $S$.

$\ul{3\Rightarrow 1}$. Similar to the proof of \leref{2.10} and \thref{4.6}.
\end{proof}

Let us now describe the connection between local comultiplications and
local multiplications. If $\Cc$ is an $A$-bimodule, and $\Delta:\ \Cc\to \Cc\ot_A \Cc$
is an associative multiplication, then $*_{\Delta}$ given by
$$(f*_{\Delta} g)(c)= g(c_{(1)}f(c_{(2)}))$$
is an associative multiplication on ${}^*\Cc:={}_A\Hom(\Cc,A)$.

\begin{lemma}\lelabel{7.5}
 If $\Delta$  and $\Delta'$ are two coassociating comultiplications on 
 an $A$-bimodule $\Cc$, then the corresponding multiplications $*_\Delta$ and $*_{\Delta'}$ associate.
\end{lemma}

\begin{proof}
If we write $\Delta(c)=c_{(1)}\otimes c_{(2)}$ and $\Delta'(c)=c_{(1')}\otimes c_{(2')}$, then the fact that $\Delta$  and $\Delta'$ coassociate means that 
$$c_{(1)}\otimes c_{(2)(1')}\otimes c_{(2)(2')}=c_{(1')(1)}\otimes c_{(1')(2)}\otimes c_{(2')}.$$
An easy computation shows that 
\begin{eqnarray*}
(f*_\Delta g)*_{\Delta'} h(c)&=&h(c_{(1)}g(c_{(2)(1')}f(c_{(2)(2')})));\\
f*_\Delta (g*_{\Delta'} h)(c)&=&h(c_{(1')(1)}g(c_{(1')(2)}f( c_{(2')}))),
\end{eqnarray*}
and the result follows.
\end{proof}

\begin{theorem}\thlabel{7.6}
If $\Cc$ is an $A$-coring with right (resp. left) weak (resp. strong) local comultiplications, then $\*C$ has right (resp. left) weak (resp. strong) local multiplications on $\Cc$.
\end{theorem}

\begin{proof}
Let $\Dd$ be a finitely generated $A$-submodule of $\Cc$ and let $\Delta$ be
a weak right local $\varepsilon$-comultiplication on $\Dd$. For $c\in\Cc$ and $f\in\*C$ we define $c\cdot_\Delta f=(\Cc\otimes_A f)\circ\Delta(c)$. An easy computation shows that this is a weak local $\varepsilon$-multiplication on $\Cc$.

Now let $\Delta$ be strong. 
If we construct $\Ee$ as in \thref{4.6}, then the projection on $\Ee$ is given exactly by right multiplication with $\varepsilon$.  
\end{proof}

\subsubsection*{Acknowledgement} The author would like to thank E. De Groot and his supervisor S. Caenepeel for helpfull comments.

\end{document}